\documentstyle[12pt,amsfonts,amssymb,leqno]{article}

\newfont\got{eufm10}

\newtheorem{proposition}{Proposition}[section]
\newtheorem{thm}[proposition]{Theorem}

\newcounter{secnum}

\begin{document}
\setcounter{section}{+1}

\begin{center}
{\Large \bf Iterations of $V$ and the core model}
\end{center}
\begin{center}
{\fbox{Preliminary version}}
\end{center}

\begin{center}
\renewcommand{\thefootnote}{\fnsymbol{footnote}}
\renewcommand{\thefootnote}{arabic{footnote}}
\renewcommand{\thefootnote}{\fnsymbol{footnote}}
{\large Ralf Schindler}
\renewcommand{\thefootnote}{arabic{footnote}}
\end{center}
\begin{center} 
{\footnotesize
{\it Institut f\"ur Formale Logik, Universit\"at Wien, 1090 Wien, Austria}} 
\end{center}

\begin{center}
{\tt rds@logic.univie.ac.at}

{\tt http://www.logic.univie.ac.at/${}^\sim$rds/}\\
\end{center}


Fix $\Omega$, a measurable cardinal. Suppose that $K^c$ (built up to $\Omega$) is
tame, and that the presence of $Q$-structures (cf. \cite{CMMW} Def. 2.1) induces an
$\Omega$ (and hence $\Omega+1$) iteration strategy for $K^c$. Then $K$ exists 
(and $\emptyset$ is excellent; cf. \cite{CMMW}
Theorem 2.7). Under these assumptions we then have the following.

\begin{thm}\label{main-result} 
Let $\pi \ \colon \ V \rightarrow M$ be an elementary embedding coming from a 
finite
coarse iteration tree on $V$ living on $V_\Omega$ such that $M$ is
transitive and ${}^\omega M \subset M$. Let $K^M = \pi(K)$ be the core model of $M$.
Then $K^M$ is an iterate of $K$, i.e., there is an iteration tree ${\cal T}$ on $K$
of successor length $\leq \Omega+1$ such that ${\cal M}^{\cal T}_\infty = K^M$.
Moreover, we'll have that
$\pi^{\cal T}_{0 \infty} = \pi \upharpoonright K$. 
\end{thm}

{\sc Proof} of \ref{main-result}. 
Let us fix $\pi \colon V \rightarrow M$ as in the statement of \ref{main-result}
throughout this proof. Let ${\cal T}$ and ${\cal T}'$ denote the iteration trees on
$K$ and $K^M$, resp., arising from the comparison of $K$ with $K^M$. 
For any ordinal $\nu$, let us say that ${\cal T}'$ is {\em beyond} $\nu$
iff for all $\alpha+1 < lh({\cal T}')$ do we have that $E_\alpha^{{\cal
T}'} \not= \emptyset \Rightarrow lh(E_\alpha^{{\cal
T}'}) \geq \nu$.

For any $\nu \in OR$, let us denote by \boldmath $(1)$\unboldmath${}_\nu$ the
claim that ${\cal T}'$ is beyond $\nu +1$.

By the argument for \cite{CMIP} Lemma 7.13, in order to show \ref{main-result} it will
suffice to prove that \boldmath $(1)$\unboldmath${}_\nu$ holds for every $\nu \in OR$.

Fix $\nu \in OR$ for a moment, and suppose that ${\cal T}'$ is beyond $\nu$
(i.e., that \boldmath $(1)$\unboldmath${}_{\bar \nu}$ holds for all ${\bar \nu} <
\nu$). 
Let $(\kappa_i \colon i<\theta)$ be the order preserving enumeration of
the set of cardinals of $K^M|\nu$,
and let $\lambda_i = \kappa_i^{+K^M|\nu}$ for each $i<\theta$. (If $\theta$ is a
successor ordinal we understand that $\lambda_{\theta-1} = \nu$.)
Let, for $i+1<\theta$, $\beta(i)$ be the least $\beta$ such that ${\cal M}_\beta^{\cal
T} | \lambda_i = K^M | \lambda_i$, and let ${\cal P}_i$ be the
longest initial segment ${\cal P}'$ of ${\cal M}_{\beta(i)}^{\cal
T}$ such that ${\cal P}(\kappa_i) \cap {\cal P}' = {\cal P}(\kappa_i) \cap
K|\nu$. (Then ${\cal P}_i$ either is a weasel, or else 
$\rho_\omega({\cal P}_i) \leq \kappa_i$.)
Let ${\vec {\cal P}} = {\vec {\cal P}}(\nu)$ denote 
the phalanx $$(({\cal P}_i \colon i+1 <
\theta)^\frown K^M,(\lambda_i \colon i+1<\theta)).$$ We let ${\vec {\cal P}}(\nu)$
be undefined if ${\cal T}'$ is not beyond $\nu$. 

For any $\nu \in OR$, let us denote by \boldmath $(2)$\unboldmath${}_\nu$ the
claim that IF ${\cal T}'$ is beyond $\nu$ THEN ${\vec {\cal P}}(\nu)$ is iterable.

Suppose now that we can prove, for any $\nu \in OR$, that $\forall {\bar \nu} < \nu \ $
\boldmath $(1)$\unboldmath${}_{\bar \nu} \ \wedge \ $ \boldmath 
$(2)$\unboldmath${}_\nu \ \Rightarrow \ $ \boldmath $(1)$\unboldmath${}_\nu$
as well as $\forall {\bar \nu} < \nu \ $
\boldmath $(2)$\unboldmath${}_{\bar \nu} \ 
\Rightarrow \ $ \boldmath $(2)$\unboldmath${}_\nu$ hold. Then 
\boldmath $(1)$\unboldmath${}_\nu$ holds for every $\nu \in OR$, and thus
\ref{main-result} is proven as pointed out above.

Now fix an ordinal $\nu$ throughout the rest of this proof. Standard
arguments\footnote{{\it NB}. We'll also have to consider a possibility discussed in
\cite{dec}.} easily
give a proof of $\forall {\bar \nu} < \nu \ $
\boldmath $(1)$\unboldmath${}_{\bar \nu} \ \wedge \ $ \boldmath 
$(2)$\unboldmath${}_\nu \ \Rightarrow \ $ \boldmath $(1)$\unboldmath${}_\nu$. We are
hence left with having to prove that $\forall {\bar \nu} < \nu \ $
\boldmath $(2)$\unboldmath${}_{\bar \nu} \ 
\Rightarrow \ $ \boldmath $(2)$\unboldmath${}_\nu$ holds; this statement can now be
rephrased as follows.

\bigskip
\boldmath $(\star) \ $\unboldmath Suppose that ${\cal T}'$ is beyond $\nu$. Suppose
further that ${\vec {\cal P}}({\bar \nu})$ is iterable for every ${\bar \nu} < \nu$.
Then ${\vec {\cal P}}(\nu)$ is iterable.

\bigskip
We are now going to prove \boldmath $(\star)$\unboldmath . Assume that the
hypotheses of \boldmath $(\star) \ $\unboldmath are met. 
We have to show that ${\vec {\cal P}}(\nu)$ is iterable. 
Let $\theta$, $(\kappa_i \colon i<\theta)$, and $(\lambda_i \colon i<\theta)$
be as in the definition of ${\vec {\cal P}}(\nu)$ given above. Let us again write 
${\vec {\cal P}}$ for ${\vec {\cal P}}(\nu)$.

Before commencing with proving anything, let us isolate three claims.
Let $\sigma \colon H \rightarrow H_{\Omega+1}$ be elementary, where
$H$ is countable and transitive, and $ran(\sigma)$
contains all the sets of current interest. Notice that $\{ H , \sigma
\upharpoonright \sigma^{-1}(K^M) \} \subset M$ by ${}^\omega M \subset M$.

Let $i+1 \in \theta \cap ran(\sigma)$. Set ${\bar {\cal P}}_i =
\sigma^{-1}({\cal P}_i)$, and let $${\cal Q}_i =
Ult({\bar {\cal P}}_i,\sigma \upharpoonright {\bar {\cal P}}_i |
\sigma^{-1}(K^M|\lambda_i))$$ be the ``lift up'' of ${\bar {\cal P}}_i$ by the
appropriate restriction of $\sigma$, which also comes with a canonical embedding
$\sigma_i \colon {\cal Q}_i \rightarrow {\cal P}_i$.
Set ${\tilde \lambda}_i = sup \ ran(\sigma \upharpoonright \sigma^{-1}(\lambda_i))$.

Let ${\vec {\cal Q}}$ denote the 
phalanx\footnote{In this sketch we simply ignore the possibility that ${\cal Q}_i$
might be a protomouse rather than a premouse. We can deal with this possibility
in the manner of \cite{covering}.}
$$(({\cal Q}_i \colon i+1 \in \theta \cap ran(\sigma))^\frown K^M,({\tilde 
\lambda}_i \colon 
i+1 \in \theta \cap ran(\sigma)).$$
It will be crucial to notice that in fact 
${\vec {\cal Q}}$ is an element of $M$.

\bigskip
{\bf Claim 1.} $((K^M,{\cal P}_i),{\tilde 
\lambda}_i)$ is iterable for each $i+1 \in \theta \cap ran(\sigma)$.

\bigskip
{\bf Claim 2.} $((K^M,{\cal Q}_i),{\tilde 
\lambda}_i)$ is iterable for each $i+1 \in \theta \cap ran(\sigma)$.

\bigskip
{\bf Claim 3.} ${\vec {\cal Q}}$ is iterable.   

\bigskip
We are now going to prove Claims 1, 2, and 3 (in that order). This will certainly
suffice as we could have thrown a potential witness to the non-iterability of ${\vec
{\cal P}}$ into $ran(\sigma)$.

\bigskip
{\sc Proof} of Claim 1. Fix $i+1 \in \theta \cap ran(\sigma)$. 
Let $j \leq i$. By our ``inductive
hypothesis,'' the phalanx $${\vec {\cal P}}(\lambda_{j+1}) 
= (({\cal P}_k \colon k \leq j)^\frown K^M,(\lambda_k \colon
k \leq j))$$ is iterable. This gives us an iterate $K^\star_j$ of $K$ together with an
embedding $\rho_j \colon K^M \rightarrow K_j^\star$ such that $\tau_j
\upharpoonright {\tilde \lambda}_j = id$. 

But now the phalanx $((K_j^\star \colon j \leq i)^\frown 
{\cal P}_i,(\lambda_j \colon
j < i)^\frown {\tilde \lambda}_i)$ is certainly iterable. 
Using the maps $\rho_j$ one can
then deduce that $((K^M,{\cal P}_i),{\tilde 
\lambda}_i)$ is iterable.

\bigskip
{\sc Proof} of Claim 2. This is a straightforward consequence of Claim 1, using the
maps $\sigma_i$.

\bigskip
{\sc Proof} of Claim 3. Let $i+1 \in \theta \cap ran(\sigma)$.
By Claim 2, $((K^M,{\cal Q}_i),{\tilde 
\lambda}_i)$ is iterable, a fact which relativizes down to $M$.
Thus by coiterating $K^M$ with $((K^M,{\cal Q}_i),{\tilde 
\lambda}_i)$ inside $M$ we get an iterate ${\cal Q}_i^\star$ of $K^M$ together with an
embedding $\tau_i \colon {\cal Q}_i \rightarrow {\cal Q}_i^\star$ such that $\tau_i
\upharpoonright {\tilde \lambda}_i = id$.

But now the phalanx $(({\cal Q}_i^\star \colon i+1 
\in \theta \cap ran(\sigma))^\frown K^M,({\tilde 
\lambda}_i \colon 
i+1 \in \theta \cap ran(\sigma)))$ is certainly iterable (in $M$, and hence in $V$).
Using the maps $\tau_i$ one can then deduce that ${\vec {\cal Q}}$ is iterable.

\bigskip
The argument given here can also be used to show that ``${\bar K}$ doesn't move in the
comparison with $K$ in the covering argument.''


\begin{thebibliography}{99}
\bibitem{covering} Mitchell, W., Schimmerling, E., and Steel, J., {\it The covering
lemma up to a Woodin cardinal}.
\bibitem{dec} Schindler, R.-D., Steel, J., and Zeman, M., {\it Deconstructing inner
model theory}, submitted.
\bibitem{CMIP} Steel, J., {\it The core model iterability problem},
Springer Verlag 1996.
\bibitem{CMMW} Steel, J., {\it Core models with more Woodin cardinals}.
\end{thebibliography}
\end{document}